\documentclass[11pt]{article}
\usepackage{amsmath}
\usepackage{amsmath,amssymb,amsfonts,amsthm,fancyhdr}
\usepackage{epsfig,graphicx,picins,picinpar,subfigure}
\usepackage{pstricks}
\usepackage{fancyvrb}
\usepackage[numbers,sort&compress]{natbib}
\begin{document}

\title{\textbf{Revisiting the  number of simple $K_4$-groups} \thanks{ $^{\dagger}$ Corresponding}}
\author{\emph{Shaohua Zhang $^1$ Wujie Shi $^2$$^{\dagger}$}}
\date{{\small  $^1$School of Mathematics and Computer Science, Yangtze Normal University, Fuling, Chongqing, 408102, China \\
E-mail address: safeprime@163.com\\
$^2$College of mathematics and Finance, Chongqing University of Arts and Sciences, Yongchuan, Chongqing, 402160, China \\
E-mail address: shiwujie@gmail.com
}}
 \maketitle

\begin{abstract}
In this paper, by solving Diophantine equations involving simple $K_4$-groups, we will
try to point out that it is not easy to prove the infinitude of simple $K_4$-groups. This problem goes far beyond what is known about Dickson's conjecture at present.

\vspace{3mm}\textbf{Keywords:} exponential Diophantine equation, simple $K_4$-group,
Dickson's conjecture, Zsigmondy's theorem

\vspace{3mm}\textbf{2000 MR  Subject Classification:}\quad 11D61,
11D45, 20D05
\end{abstract}


\section{Introduction}
\setcounter{section}{1}\setcounter{equation}{0}

In the famous book \emph{Unsolved Problems in Group Theory}, the following problem is asked: is the number of simple $K_4$-groups finite or infinite? See \cite {1}: Problem 13.65. This problem is the first to be posed by the second author W. Shi in 1991 \cite {2}. Denote by $\mathbb{N}$ and $\mathbb{P}$
the sets of positive integers and the set of prime numbers, respectively. In \cite {2}, the second author claimed that the simple $K_4$-group problem can be reduced to
the four Diophantine problems:\\
(1) $p^2-1=2^a3^bq^c$, $p, q\in \mathbb{P}$, $p>3$, $q>3$, $a,b,c\in \mathbb{N}$, \\
(2) $2^m-1=p$, $2^m+1=3q^n$, $p, q\in \mathbb{P}$, $p>3$, $q>3$, $m,n\in \mathbb{N}$, \\
(3) $3^m-1=2p^n$, $3^m+1=4q$, $p, q\in \mathbb{P}$, $p>3$, $q>3$, $m,n\in \mathbb{N}$, \\
(4) $3^m-1=2p$, $3^m+1=4q^n$, $p, q\in \mathbb{P}$, $p>3$, $q>3$, $m,n\in \mathbb{N}$.

\vspace{3mm}In 2001, Y. Bugeaud, Z. Cao and M. Mignotte \cite {3} showed that if $n>1$, (2) and (4) have no solution and (3) has only the solution $(p,q,m,n)=(11,61,5,2)$.

\vspace{3mm}In Section 2, we will prove that if $c>1$, (1) has only the solutions $(p,q,a,b,c)=(97,7,6,1,2)$ and $(p,q,a,b,c)=(577,17,7,2,2)$. Our methods are slightly different from those in \cite {3}.

\vspace{3mm}Thus, the infinitude of simple $K_4$-groups can be decided by
the following three Diophantine problems:\\
(5) $p^2-1=2^a3^bq$, $p, q\in \mathbb{P}$, $p>3$, $q>3$, $a,b\in \mathbb{N}$, \\
(6) $2^m-1=p$, $2^m+1=3q$, $p, q\in \mathbb{P}$, $p>3$, $q>3$, $m\in \mathbb{N}$, \\
(7) $3^m-1=2p$, $3^m+1=4q$, $p, q\in \mathbb{P}$, $p>3$, $q>3$, $m\in \mathbb{N}$.

\vspace{3mm}\noindent{\bf Remark:~~}%
In \cite {3}, one would see that the infinitude of simple $K_4$-groups can be decided by
the  Diophantine equation: $q(q^2-1)=\gcd (2,q-1)2^a3^bp^cr^d$,  where $q$ is a prime power and $3<p,r\in \mathbb{P}$, $r\neq p$, $a,b,c,d\in \mathbb{N}$. However, (5), (6), (7) and this result are consistent on the infinitude of simple $K_4$-groups.

\vspace{3mm}In this paper, by considering Diophantine equations (5), (6) and (7), our main aim is to try to point out that it is very difficult to determine the infinitude of simple $K_4$-groups, and this problem goes far beyond what is known about Dickson's conjecture \cite {4} at present. For the details, see  Section 3.

\section{Main theorem and its Proof}
\vspace{3mm}\noindent{\bf  Theorem 1:~~}%
If $1<c\in \mathbb{N}$, Diophantine equation $p^2-1=2^a3^bq^c$ has only the solutions $(p,q,a,b,c)=(97,7,6,1,2)$ and $(p,q,a,b,c)=(577,17,7,2,2)$ satisfying $p, q\in \mathbb{P}$, $p>3$, $q>3$, $a,b\in \mathbb{N}$.

\vspace{3mm}\noindent{\bf  Lemma 1 (Zsigmondy's theorem  \cite {5}):~~}%
If $a>b>0$, $\gcd (a,b)=1$ and $n>1$ are positive integers, then $a^n+b^n$ has a prime factor that does not divide $a^k+b^k$ for all positive integers $k<n$, with the exception $2^3+1^3$;  $a^n-b^n$ has a prime factor that does not divide $a^k-b^k$ for all positive integers $k<n$ unless $a=2,b=1$ and $n=6$; or $a+b$ is a power of $2$ and $n=2$.

\vspace{3mm}\noindent{\bf  Lemma 2 (\cite {6}):~~}%
The Diophantine equation $3^m-2y^q=1$ has only the solution $(y,m,q)=(11,5,2)$ satisfying $2<m\in \mathbb{N}$, $2\leq y,q \in \mathbb{N}$.

\vspace{3mm}\noindent{\bf  Lemma 3 (\cite {7}):~~}%
The Diophantine equation $x^m-y^n=1$ has only the solution $(x,y,m,n)=(3,2,2,3)$ satisfying $1<m,n\in \mathbb{N}$, $x,y \in \mathbb{P}$.

\vspace{3mm}\noindent{\bf  Lemma 4 (\cite {8}):~~}%
The Diophantine equation $x^2+1=2y^n$ has only the solution $(x,y,n)=(239,13,4)$ satisfying $x,y,n\in \mathbb{N}$, $y>1,n>2$.

\vspace{3mm}\noindent{\bf  Lemma 5 (\cite {9}):~~}%
The Diophantine equation $3x^2+1=4y^n$ with $1<n\equiv 1 (\mod2)$ and $x,y,n\in \mathbb{N}$
has only the solution $x=y=1$.

\vspace{3mm}\noindent{\bf  Lemma 6 (\cite {3}):~~}%
The Diophantine equation $2^m+1=3y^q$ with $1<m,1<y,1<q\in \mathbb{N}$
has no solution.

\vspace{3mm}\noindent{\bf  Proof of theorem 1:~~} %
If $1<c\in \mathbb{N}$, Diophantine equation $p^2-1=2^a3^bq^c$ has solutions $(p,q,a,b,c)$ such that $p, q\in \mathbb{P}$, $p>3$, $q>3$, $a,b\in \mathbb{N}$, then $a\geq 3$ and $p^2\equiv 1 (\mod 2^a)$. Notice that if $a\geq 3$, $x^2\equiv 1 (\mod 2^a)$ has only four solutions, say $x=\pm 1, \pm 1+2^{a-1}$. Therefore, we must have $p=\pm 1 +2^{a-1}+k 2^a$ with $k \in \mathbb{N}\cup \{0\}$. When $k \in \mathbb{N}$,  $p^2-1=2^a3^bq^c$ can be reduced to
the following Diophantine problems with $3\leq a,1\leq b,1<c\in \mathbb{N}$, $p, q\in \mathbb{P}$, $p>3$, $q>3$:\\
(8) $2k+1=3^b$, $ q^c=1 +2^{a-2}+k 2^{a-1}$, $p=1 +2^{a-1}+k 2^a$,  \\
(9) $2k+1=q^c$, $ 3^b=1 +2^{a-2}+k 2^{a-1}$, $p=1 +2^{a-1}+k 2^a$,   \\
(10) $2k+1=3^b$, $ q^c=-1 +2^{a-2}+k 2^{a-1}$, $p=-1 +2^{a-1}+k 2^a$, \\
(11) $2k+1=q^c$, $ 3^b=-1 +2^{a-2}+k 2^{a-1}$, $p=-1 +2^{a-1}+k 2^a$.

\vspace{3mm}Rewrite these equations, we get (with $3\leq a,1\leq b,1<c\in \mathbb{N}$, $p, q\in \mathbb{P}$, $p>3$, $q>3$):\\
(12) $ q^c-1 =2^{a-2}3^b$, $p=-1 +2q^c$,  \\
(13) $ 3^b-1= 2^{a-2}q^c$, $p=-1 +2\times3^b$,  \\
(14) $ q^c +1 =2^{a-2}3^b$, $p=1 +2q^c$,  \\
(15) $ 3^b+1= 2^{a-2}q^c$, $p=1 +2\times3^b$.

\vspace{3mm}Next, we will prove that if $1<c$, then (13), (14), (15) have no solution and (12) has only two solutions satisfying the conditions.

\vspace{3mm} Clearly, if (14) has solutions, then $c\neq 2$. By Lemma 1, $q^c +1$ has at least one prime factor $m$ that does not divide $q^r +1$ for all positive integers $r<c$. However, $m\neq 2,3$. This leads to a contradiction since $ q^c +1 =2^{a-2}3^b$.

\vspace{3mm} Now, let's consider (15). Assume that $ 3^b+1= 2^{a-2}q^c$ has solutions. If $b$ is even, then $a=3$. By Lemma 4, it is impossible. Hence, $b$ is odd. We deduce that $a=4$ and get that $ 3^b+1= 4q^c$. If $c$ is even, then $ 3^b+1= 4q^c$ has no solution with $c>1$ (\cite {10}). If $c$ is odd, then $ 3^b+1= 4q^c$ has no solution by Lemma 5 (since $c>1$). So, (15) has no solution.

\vspace{3mm} Suppose that (13) $ 3^b-1= 2^{a-2}q^c$ has solutions. If $b$ is odd, then $a=3$ and $ 3^b-1= 2q^c$ has no solution such that $p=-1 +2\times3^b$ by Lemma 2. Let $b$ be even. Write $b=2r$. We obtain that $a\geq 5$ and $\frac{3^r-1}{2}\frac{3^r+1}{2}=2^{a-4}q^c$. By Lemma 3, one can prove that $\frac{3^r-1}{2}\frac{3^r+1}{2}=2^{a-4}q^c$ has no solution. Thus, (13) has no solution.

\vspace{3mm} If $ q^c-1 =2^{a-2}3^b$ has solutions with $c>1$, then $c=2$  by Lemma 1. By Lemma 3, one can obtain that (12) has only solutions $(q,a,b)=(7,6,1)$ and $(q,a,b)=(17,7,1)$.  It leads that  Diophantine equation $p^2-1=2^a3^bq^c$ has only solutions $(p,q,a,b,c)=(97,7,6,1,2)$ and $(p,q,a,b,c)=(577,17,7,2,2)$ satisfying $p, q\in \mathbb{P}$, $p>3$, $q>3$, $a,b,c\in \mathbb{N}$ and $c>1$.

\vspace{3mm}Finally, we consider the case $k=0$. Obviously, $p^2-1=2^a3^bq^c$ can be reduced to
the following Diophantine problems with $3\leq a,1\leq b,1<c\in \mathbb{N}$, $p, q\in \mathbb{P}$, $p>3$, $q>3$:\\
(16) $ 3^bq^c=1 +2^{a-2}$, $p=1 +2^{a-1}$,  \\
(17) $ 3^bq^c=-1 +2^{a-2}$, $p=-1 +2^{a-1}$.

\vspace{3mm} By Lemma 3, if (16) or (17) has solutions, then $b$ must be $1$. Furthermore, using Lemma 6, one can show that (16) and (17) have no solution satisfying the conditions. This proves Theorem 1.

\section{Our conclusion}
In this section, we will
try to point out that it is not easy to prove the infinitude of simple $K_4$-groups. By the aforementioned results, one will see that  the infinitude of simple $K_4$-groups can be decided by
the following  Diophantine problems:\\
(18) $2^m-1=p$, $2^m+1=3q$, $p, q\in \mathbb{P}$, $p>3$, $q>3$, $m\in \mathbb{N}$, \\
(19) $3^m-1=2p$, $3^m+1=4q$, $p, q\in \mathbb{P}$, $p>3$, $q>3$, $m\in \mathbb{N}$, \\
(20) $ q-1 =2^{a-2}3^b$, $p=-1 +2q$, $p, q\in \mathbb{P}$, $a\geq 3$, $a,b\in \mathbb{N}$, \\
(21) $ 3^b-1= 2^{a-2}q$, $p=-1 +2\times3^b$,   $p, q\in \mathbb{P}$, $a\geq 3$, $a,b\in \mathbb{N}$, \\
(22) $ q +1 =2^{a-2}3^b$, $p=1 +2q$,   $p, q\in \mathbb{P}$, $a\geq 3$, $a,b\in \mathbb{N}$, \\
(23) $ 3^b+1= 2^{a-2}q$, $p=1 +2\times3^b$,  $p, q\in \mathbb{P}$, $a\geq 3$, $a,b\in \mathbb{N}$,\\
(24) $ 3q=1 +2^{a-2}$, $p=1 +2^{a-1}$,  $p, q\in \mathbb{P}$, $a\geq 3$, $a\in \mathbb{N}$,\\
(25) $ 3q=-1 +2^{a-2}$, $p=-1 +2^{a-1}$.  $p, q\in \mathbb{P}$, $a\geq 3$, $a\in \mathbb{N}$.

\vspace{3mm} Hence, if the number of simple $K_4$-groups is infinite, then one of the following holds, which goes far beyond what is known about Dickson's conjecture at present:\\
(26) $f_1(x)=x$ and $f_2(x)=3x-2$ represent simultaneously primes for infinitely many integers $x$ by (18), \\
(27) $f_1(x)=x$ and $f_2(x)=2x-1$ represent simultaneously primes for infinitely many integers $x$ by (19) or (20),\\
(28) $f_1(x)=x$ and $f_2(x)=2x+1$ represent simultaneously primes for infinitely many integers $x$ by (22),\\
(29) $f_1(x)=x$ and $f_2(x)=4x+1$ represent simultaneously primes for infinitely many integers $x$ by (21) (Note that  by (21) one can deduce that $a$ must be $3$.),\\
(30) $f_1(x)=x$ and $f_2(x)=2^{a-1}x-1$ represent simultaneously primes for infinitely many integers $x$ by (23), where $a=3$ or $a=4$,\\
(31) $f_1(x)=x$ and $f_2(x)=6x-1$ represent simultaneously primes for infinitely many integers $x$ by (24), \\
(32) $f_1(x)=x$ and $f_2(x)=6x+1$ represent simultaneously primes for infinitely many integers $x$ by (25).

\vspace{3mm}\noindent{\bf  Dickson��s conjecture:~~}%
 Let $1\leq s \in \mathbb{N}$, $f_i(x) = a_i +b_ix$ with $a_i$ and $b_i$ integers,
$b_i \geq 1$ (for $i = 1, ..., s $). If there does not exist any integer $n > 1$ dividing
all the products $\prod _{i=1}^{i=s}f_i(k)$, for every integer $k$, then there exist infinitely
many natural numbers $m$ such that all numbers $f_1(m), ..., f_s(m)$ are primes.

\vspace{3mm}The case $s = 1$ is Dirichlet's theorem. Two  special cases  are well known conjectures: there are infinitely many twin primes ( $f_1(x)=x$ and $f_2(x)=x+2$ represent simultaneously primes for infinitely many integers $x$), and there are infinitely many Sophie Germain primes ($f_1(x)=x$ and $f_2(x)=2x+1$ represent simultaneously primes for infinitely many integers $x$). Unfortunately, as we know, even these two simple cases, nobody has proved up until now, let alone the following special cases: (26), (27), (29), (30), (31) and (32).

\vspace{3mm}Dickson's conjecture is further extended by Schinzel's hypothesis H \cite {11}. This research line is to try to generalize Dickson��s conjecture (or Schinzel's hypothesis H) to the case of multivariable polynomials and obtain a similarly quantitative form of Schinzel-Sierpinski's
Conjecture---Bateman-Horn's conjecture\cite {12}.

In  \cite {13}, the first author guesses that if the linear polynomial map $$F(x)=\left\{
\begin{array}{c}
f_1(x_1,...,x_n)=a_{11}x_1+...+a_{1n}x_n +b_1\\
...........................................................\\
f_m(x_1,...,x_n)=a_{m1}x_1+...+a_{mn}x_n +b_m\\
\end{array}
\right.$$ on $\mathbb{Z}^n$ is
admissible, then polynomials
$f_1(x_1,...,x_n),...,f_m(x_1,...,x_n)$  with integral coefficients represent simultaneously
prime numbers for infinitely many integral points $(x_1,...,x_n)$, where $\mathbb{Z}$ is the set of all integers.

\vspace{3mm}\noindent{\bf Remark:~~}%
We call the linear polynomial map $$F(x)=\left\{
\begin{array}{c}
f_1(x_1,...,x_n)=a_{11}x_1+...+a_{1n}x_n +b_1\\
...........................................................\\
f_m(x_1,...,x_n)=a_{m1}x_1+...+a_{mn}x_n +b_m\\
\end{array}
\right.$$  on $\mathbb{Z}^n$ is
admissible if for every positive integer $r$ there exists an
integral point $x=(x_1,...,x_n)\in Z^n$  such that
$r$ is coprime to $f_1(x)\times...\times
f_m(x)=\prod_{i=1}^{i=m}f_i(x)$, moreover, $f_i(x)>1$ for $1\leq
i\leq m$.

In 2006 \cite {14}, B. Green and T. Tao  generalized Hardy-Littlewood Conjecture. Roughly speaking:

let $F(x)=\left\{
\begin{array}{c}
f_1(x_1,...,x_n)=a_{11}x_1+...+a_{1n}x_n +b_1\\
...........................................................\\
f_m(x_1,...,x_n)=a_{m1}x_1+...+a_{mn}x_n +b_m\\
\end{array}
\right.$ on $\mathbb{Z}^n$ be the linear polynomial map, if  polynomials
$f_1(x_1,...,x_n),...,f_m(x_1,...,x_n)$  with integral coefficients represent simultaneously
prime numbers for infinitely many integral points $(x_1,...,x_n)$, there must be a constant $C_F$ only depending on $F(x)$, such that when $h\rightarrow\infty$, $$\# \{x\in [-h,h]^n|F(x)\in \mathbb{P}^m\} \approx C_F\frac{h^n}{\log^mh}.$$

This formula implies Prime Number Theorem, Green-Tao Theorem \cite {15} which states that
$$F(x)=\left\{
\begin{array}{c}
f_1(x_1,x_2)=x_1\\
f_2(x_1,x_2)=x_1+x_2\\
.........................\\
f_m(x_1,x_2)=x_1+(m-1)x_2\\
\end{array}
\right.$$  represent simultaneously
prime numbers for infinitely many integral points $(x_1,x_2)$, moreover, there must be a constant $C_F$ only depending on $F(x)$, such that when $h\rightarrow\infty$, $$\# \{x\in [-h,h]^2|F(x)\in \mathbb{P}^m\} \approx C_F\frac{h^2}{\log^mh},$$
 and so on.

 If the aforementioned formula $\# \{x\in [-h,h]^n|F(x)\in \mathbb{P}^m\} \approx C_F\frac{h^n}{\log^mh}$ is true, maybe, people would call it Prime Point Theorem. From Prime Number Theorem to  Prime Point Theorem, it will be a very long way.

It is worthwhile to point out a latest advance of Dickson's conjecture. The new result in 2013, from Yitang Zhang of the University of New Hampshire in Durham, finds that there are infinitely many pairs of primes that are less than 70 million without relying on unproven conjectures \cite {16}. This implies that in the following 35000000 systems of equations:

 $\left\{
\begin{array}{c}
f_1(x)=x   \\
g_1(x)=x+2\\
\end{array}
\right.$, $\left\{
\begin{array}{c}
f_2(x)=x   \\
g_2(x)=x+4\\
\end{array}
\right.$,......, $\left\{
\begin{array}{c}
f_{35000000 }(x)=x   \\
g_{35000000 }(x)=x+70000000 \\
\end{array}
\right.$

 there is at least a system of equations, say $\left\{
\begin{array}{c}
f_i(x)=x   \\
g_i(x)=x+2i\\
\end{array}
\right.$, such that $f_i(x)$ and $g_i(x)$ represent simultaneously
prime numbers for infinitely many integer numbers $x$.

On the other hand, even if Dickson's conjecture holds, it is not obvious that the number of simple $K_4$-groups is infinite. In fact, the number of simple $K_4$-groups is very closely related to the  following problem: let $$F(x)=\left\{
\begin{array}{c}
f_1(x)=a_{1}x +b_1\\
....................\\
f_m(x)=a_{m}x+b_m\\
\end{array}
\right.$$ on an infinite subset $A$ of $\mathbb{Z}$ be a linear polynomial map, what
condition does $F(x)$ satisfy so that $f_1(x),...,f_m(x)$ represent simultaneously
prime numbers for infinitely many integer numbers $x$? Unfortunately, on this problem, we cannot find any references. It Looks harder than Dickson's conjecture. Anyway, due to that fact it is  closely tied with many topics in number theory such as Fermat's primes, Mersenne primes, Dickson's conjecture and so on, we think that  determining the number of simple $K_4$-groups is significative. All these and related questions, we hope to further investigate.

\section{Acknowledgements}
This work was partially supported by the NSFC(No.11171364), the Scientific Research Foundation of Yangtze Normal University, the Science and Technology Research Project of Chongqing City Board of  Education
(No.KJ121316), the Doctoral Fund of Ministry of Education of China (No.20090131120012) and the Innovation Fundation of Chongqing(KJTD201321).

\clearpage
\end{document}